\documentclass[runningheads]{llncs}

\usepackage{graphicx}
\usepackage{multirow}
\usepackage{wrapfig}
\usepackage{amsmath, amssymb, amsthm}
\usepackage{booktabs}

\usepackage{hyperref}

\hypersetup{
    colorlinks=TRUE,
    citecolor=gray,
    pdftitle ={proof},
    linkcolor=gray,
    filecolor=gray,      
    urlcolor=gray
}

\usepackage[british]{babel} 
\usepackage[color=lightgray]{todonotes}

\begin{document}
\title{Truth, Proof, and Reproducibility: There's no counter-attack for the codeless}
%
%
\author{Charles T. Gray\inst{1}\orcidID{0000-0002-9978-011X}\thanks{Thank you to Kerrie Mengersen, Kate Smith-Miles, Mark Padgham, Hien Nguyen, Emily Kothe, Fiona Fidler, Mathew Ling, Luke Prendergast, Adam Sparks, Hannah Fraser, Felix SingletonThorn, James Goldie, Michel Penguin (Michael Sumner), in no particular order, with whom initial bits and pieces of this paper were discussed. Special thanks to 
Brian A. Davey for proofing the proofs and 
Alex Hayes for his edifying post~\cite{hayes_testing_2019}} \and
Ben Marwick\inst{2}\orcidID{0000-0001-7879-4531}}
\institute{La Trobe University, Melbourne \email{charlestigray@gmail.com}  \and
University of Washington, Seattle \email{bmarwick@uw.edu}
}
\maketitle
\begin{abstract}

Current concerns about reproducibility in many research communities can be traced back to a high value placed on empirical reproducibility of the physical details of scientific experiments and observations. For example, the detailed descriptions by 17th century scientist Robert Boyle of his vacuum pump experiments are often held to be the ideal of reproducibility as a cornerstone of scientific practice. Victoria Stodden has claimed that the computer is an analog for Boyle's pump -- another kind of scientific instrument that needs detailed descriptions of how it generates results. In the place of Boyle's hand-written notes, we now expect code in open source programming languages to be available to enable others to reproduce and extend computational experiments. In this paper we show that there is another genealogy for reproducibility, starting at least from Euclid, in the production of proofs in mathematics. Proofs have a distinctive quality of being necessarily reproducible, and are the cornerstone of mathematical science. However, the task of the modern mathematical scientist has drifted from that of blackboard rhetorician, where the craft of proof reigned, to a scientific workflow that now more closely resembles that of an experimental scientist. So, what is proof in modern mathematics? And, if proof is unattainable in other fields, what is due scientific diligence in a computational experimental environment? How do we measure truth in the context of uncertainty? Adopting a manner of Lakatosian conversant conjecture between two mathematicians, we examine how proof informs our practice of computational statistical inquiry. We propose that a reorientation of mathematical science is necessary so that its reproducibility can be readily assessed. 

\end{abstract}

In David Auburn's Pulitzer prize-winning 2000 play \emph{Proof}, a young mathematician, Catherine, struggles to prove to another mathematician, Hal, that her argument is not a reproduction of the intellectual work of her deceased father, a professor~\cite{auburn_proof_2001}. Her handwriting similar to her father's, there is no way to discern her proof from his. But if Catherine were a computational scientist, we would have a very different story. We reimagine Hal challenging Catherine for different mathematical questions and the reproducibility of her solutions. We consider simple to complex mathematical questions that can be answered at the blackboard, and then consider the scenario where Catherine must use a combination of mathematical and computational tools to answer a question in mathematical science. Via these scenarios, we question to what extent proof methodology continues to inform our choices as mathematical scientists become as much research software engineers as they are mathematicians. 

Mathematical science is the compendium of research that binds the Catherine's methodology of work indistinguishably from her father's. However, in computational science, we not only do not have a common language in the traditional sense, with programming languages such as Python, R, and C++ performing overlapping tasks, but our research workflows comprise tools and platforms and operating systems, such as Linux or Windows, as well. Many inadvertent reasons conspire so that scientists are arriving at similar problems with different approaches to data management and version control. Code scripts, arguably the most immediately analogous to mathematical proof, are but one of the many components that make up the outputs of computational science. 

If Catherine were a contemporary computational mathematician, she would not only struggle to reproduce another person's work, but she would likely struggle to reproduce her own. She may be overwhelmed by the diversity of research outputs~\cite{bryan_excuse_2018}, and find that she needs to rewrite her work to unpick what she did with specific computational functions under specific software package releases. The language of mathematical science has changed from something we write, to something we collect. In order to diligently answer scientific questions computationally, the mathematician must now consider her work within that of a research compendium. In this paper we ask: how can we extend the certainty afforded by a mathematical proof further down the research workflow into the `mangle of practice' ~\cite{pickering2010mangle}? We show that communities of researchers in many scientific disciplines have converged on a toolkit that borrows heavily from software engineering to robustly provides many points to verify certainty, from transparency via version control, to stress testing of algorithms. We focus on unit testing as a strong measure of certainty. 

\section{The technological shift in mathematical inquiry}

The task of a mathematical scientist in the pre-computer age was largely that of a blackboard rhetorician, where the craft of proof reigned. For a proof such as that featured in Auburn's play, the argument can often be included in the article, or as a supplementary file. This allows the reader to fully reproduce the author's reasoning, by tracing the flow of argument through the notation. As computers have become ubiquitous in research, mathematical scientists have seen their workflow shift to one that now more closely resembles that of a generic scientist, concerned with diligent analysis of observational and experimental data, mediated by computers~\cite{peng_reproducible_2011-1}. But the answer to the question of what constitutes a diligent attempt to answer a scientific question examined in a computationally intensive analysis, is unclear, and remains defined by the era of the blackboard mathematician.  

So, what is proof in mathematics, when experimental and computer-assisted methods are common? And, beyond mathematics, in fields where literal proofs are unattainable, what counts as an equivalent form of scientific certainty in a computational experimental environment? How do we measure truth in the context of uncertainty? Among the histories of science we can trace three efforts to tackle these questions. First is the empirical effort, most prominently represented by Robert Boyle (1627-1691), known for his vacuum pump experiments ~\cite{shapin2011leviathan}. Boyle documented his experiments in such detail and to an extent that was uncommon at the time. He was motivated by a rejection of the secrecy common in science at his time, and by a belief in the importance of written communication of experimental expertise (as a supplement to direct witnessing of experimental procedures). Boyle's distinctive approach of extensive documentation is often cited by modern advocates of computational reproducibility ~\cite{stodden_what_2014}. Making computer code openly available to the research community is argued to be the modern equivalent of Boyle's exhaustive reporting of his equipment, materials, and procedures~\cite{leveque_reproducible_2012}.   

A second effort to firming up certainty in scientific work, concerned with statistical integrity, can be traced at least as far back as Charles Babbage (1791-1871), mathematician and inventor of some of the first mechanical computers. In his 1830 book `Reflections on the Decline of Science in England, and on Some of Its Causes' he criticised some of his contemporaries, characterising them as `trimmers' and `cooks'~\cite{haack_defending_2011}. Trimmers, he wrote, were guilty of smoothing of irregularities to make the data look extremely accurate and precise. Cooks retained only those results that fit their theory and discarded the rest~\cite{merton_social_1996}. These practices are now called data-dredging, or p-hacking, where data are manipulated or removed from an analysis until a desirable effect or p-value is obtained~\cite{head_extent_2015}.

A third effort follows the history of formal logic through to the time when an equivalence between philosophical logic and computation was noted. This observation is called the Curry-Howard isomorphism or the proofs-as-programs interpretation. First stated in 1959, this correspondence proposed that proofs in some areas of mathematics, such as type theory, are exactly programs from a particular programming language ~\cite{sorensen2006lectures}. The bridging concepts come from intuitionistic logic and typed lambda calculi, which have lead to the design of computational formal proof management systems such as the Coq language. This language is designed to write mathematical definitions, execute algorithms and theorems, and check proofs ~\cite{bertot2008short}. This correspondence has not been extensively discussed in the context of reproducibility, but we believe it has relevance and is motivating beyond mathematics. Our view is that this logic-programming correspondence can be extended in a relaxed way beyond mathematics in proofs to scientific claims in general, such that computational languages can express those claims in ways that can establish a high degree of certainty.  

Questions of confidence in scientific results are far from restricted to the domains of mathematics or computers; indeed, science is undergoing a broad reexamination under what is categorised as a crisis of inference~\cite{fidler_reproducibility_2018}. How we reproduce scientific results is being examined across a range of disciplines~\cite{camerer_evaluating_2016,wallach_reproducible_2018}. An early answer to some of these questions is that authors should make available the code that generated the results in their paper~\cite{10.1093/biostatistics/kxq028,stodden_setting_2013}. These recommendations mark the emergence of a concern for computational reproducibility in mathematics. This paper extends this argument for computational reproducibility further into the workflow of modern statistical inquiry, expanding and drawing on solutions proposed by methods that privilege computational reproducibility.  

 
Systemic problems are now being recognised in the practice of conventional applied statistics, with a tendency towards \emph{dichotomania}~\cite{amrhein_scientists_2019} that reduces complex and nuanced questions to Boolean statements of \verb|TRUE| or \verb|FALSE|. This has diluted the trust the can be placed in scientific results, and led to a crisis of replication, where results can not easily be reproduced~\cite{fidler_reproducibility_2018}.

As the conventions of statistics are called into question, it stands to reason that the research practices of the discipline of statistics itself require examination. For those practicing statistical computing, a conversation is emerging about what constitutes best practice~\cite{wilson_best_2014}. But best practice may be unrealistic, especially for those applying statistics from fields where their background has afforded limited computational training. And thus the question is becoming reframed in terms of  \emph{good enough} standards~\cite{wilson_good_2017} we can reasonably request of statistical practitioners. By extension, we must reconsider how we prepare students in data-analytic degree programs. 

\begin{figure}
\centering
\includegraphics[width=\textwidth]{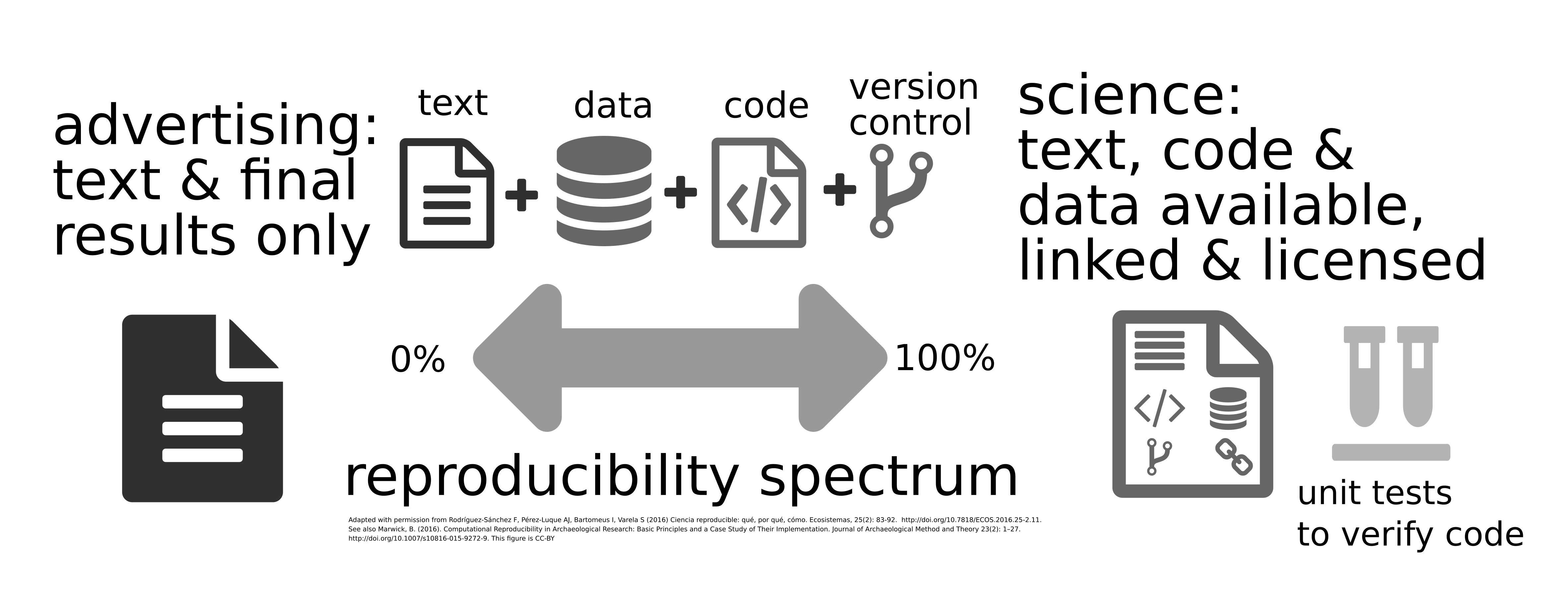}
\caption{We propose updating this spectrum of reproducibility~\cite{marwick_packaging_2018} with unit tests for data analysis. In addition to the advertising, the \textbf{formal} scientific argument put forward,  many \emph{informal} and traditionally hidden scientific outputs comprise the compendium of research that produces the results. Given the underutilised nature of unit tests, we suggest there is further work to be done to facilitate the adoption of \emph{good enough}~\cite{wilson_good_2017} research software engineering practices for answering mathematical questions computationally. The informal components of mathematical research compendium are shaded grey.} \label{fig:repro-spectrum}
\end{figure}

 Proofs, derivations, verification, all form the work of mathematics. How do we make mathematical arguments in a computational\footnote{We focus in this manuscript on R packages, but the reader is invited to consider these as examples rather than definitive guidance. The same arguments hold for other languages, such as Python, and associated tools.} environment?  In constructing mathematical arguments, we posit that we require an additional core element: unit testing for data analysis. We propose an expansion of the spectrum of reproducibility, Figure \ref{fig:repro-spectrum}, to include unit testing for data analytic algorithms facilitated by a tool such as \verb|testthat::|~\cite{wickham_testthat_2011}, for answering mathematical research questions computationally. In order to motivate this practice, we turn to the purest of sciences, mathematical proof.

\section{Truth in mathematics}

The titular proof~\cite{auburn_proof_2001} of Auburn's play is a mathematical argument, a formalised essay in mathematical science. The creator of the proof, Catherine, is questioned by Hal, who is capable of following the argument; that is, Hal can \emph{replicate} an approximation of the type of thought process that leads to a \emph{reproduction} of the argument presented in the proof. 

In Figure \ref{fig:repro-spectrum}, we have coloured the  components, black \textbf{formal} argument, and grey \emph{informal} work, of mathematics Hal would need to reproduce the proof. In order to verify the results, Hal would need to follow the formal argument, to understand what was written in the proof, but also need to do informal work, to understand the links between concepts for verification. 

Hal would come to the problem with a different background and education to Catherine. Although work is necessary for the verification of the results, the reproduction of the reasoning, the work required would be different for Hal and Catherine, based on their respective relevant preparation. However, the language of mathematics carries enough uniformity that Hal can fill in the work he requires to understand the result, from reasoning and mathematical texts. If Catherine were asking a mathematical question computationally, the presentation of the results carries not millennia of development of methodology the noble craft of mathematics, but less than a century of frequently disconnected developments separated by disparate disciplines.   

We begin with traditional mathematics and end with answering questions in computational mathematics. To this aim, we adopt, in the manner of Lakatos' \emph{Proofs and Refutations}' conversant conjecture, scenarios between Hal and Catherine, where Hal challenges Catherine over her authorship of the proof. In each scenario, we imagine the challenge would play out for different ways of answering mathematical questions. We argue the thinking work of mathematical science is not as immediately inferable in a computational experimental environment, and that the roots of mathematical science in proof lead to an overconfidence that science is as readily reproducible as a proof.    

\subsection{Prove it!}\label{sec:prove-it}

Let us suppose Catherine claimed she could demonstrate a property about the order\footnote{Let $P$ be a set. An \emph{order} on $P$ is a binary relation $\leqslant$ on $P$ such that, for all $x, y, z \in P$: we have $x \leqslant x$; with $x \leqslant y$ and $y \leqslant x$ imply $x = y$; and, finally, $x \leqslant y$ and $y \leqslant z$ imply $x \leqslant z$. We then say $\leqslant$ is reflexive, antisymmetric, and transitive, for each of these properties, respectively~\cite{davey_introduction_2002-1}.} on natural numbers, $\mathbb N =  \{1, 2, 3, \dots \}$, the counting numbers. 

The order on a set of numbers is dense if, for any two numbers we can find a number in between. More formally, we say an ordered set $P$ is dense if, for all $x < y$ in $P$, there exists $z$ in $P$ such that $x < z < y$. 

Catherine presents the following argument that the order on $\mathbb N$ is not dense. In this case she chooses a type of \emph{indirect} proof, an \emph{existence} proof~\cite{brown2013partial}, where she presents a counterexample demonstrating that the density property is not true for all cases for $\mathbb N$.   

\begin{proof}{\textbf{The order on $\mathbb N$ is not dense.}}
Let us, in the spirit of Lewis Carroll\footnote{Lewis Carroll, author of \emph{Alice in Wonderland}, is a writing pseudonym used by Charles Lutwidge Dogson, born in 1832, who taught mathematics at Christ Church, Oxford~\cite{carroll_annotated_1999}.  
}, be contrary and suppose, by way of contradiction, that the order on $\mathbb N$, \emph{is} dense. Then for any two numbers, $x$ and $y$, in $\mathbb N$ such that $x < y$, I should be able to find a distinct number $z$  in $\mathbb N$ between them, that is $x < z < y$. But, consider the numbers 3 and 4. Let $x = 3$ and $y = 4$, then $x < y$. There is no distinct number, $z$, that exists between $x$ and $y$. Since this rule must be true for any two numbers $x < y$ in the order to be dense, we have shown the order on natural numbers $\mathbb N$ is not dense.\qed
\end{proof}

A standard way to prove something is \emph{not} true, is to assume it \emph{is} true, and derive a contradiction ~\cite{davey_when_2009}. Arguably, this reasoning  goes to the heart of the problem of \emph{dichotomania} lamented by 800 scientists in a recent protest paper about the misinterpretation of statistics in \emph{Nature}~\cite{amrhein_scientists_2019}. A null hypothesis test of a difference between two groups will assume the opposite of what we suspect is true; we believe there to be a difference between two groups and take a sample from each of the groups and perform a test. This test assumes there is no difference, null, between the two groups and that any observed differences in sampling are due to random chance. The calculation returned, the $p$-value, is the likelihood we would observe the difference under those null assumptions. Crucially, the calculation returned is probabilistic, a number between 0 and 1, not a \verb|TRUE| or \verb|FALSE|, the logic of a proof by contradiction. The logic does not apply to a situation where, within a single group of people, some people might be resistant to treatment, and some might not be, say, and we have estimated a likelihood of the efficacy of the treatment. Dichotomania is the common misinterpretation of a probabilistic response in a dichotomous framework; scientists are unwittingly framing null hypothesis significance testing in terms of a proof by contradiction. 

In order to illustrate our central point, we now turn to a direct argument, rather than the indirect approach of contradiction, in order to examine the process of the making of a proof. In both the case of the direct, and indirect proofs, however, Hal could challenge Catherine, as he did in the play.
 
 \begin{quote}
     `Your dad might have written it and explained it to you later. I'm not saying he did, I'm just saying there's no proof that you wrote this'~\cite{auburn_proof_2001}.
 \end{quote}
   
\subsection{The steps in the making of a proof}

Let us now suppose Catherine's proof instead demonstrated a density property on the order on real numbers, $\dots, -3, \dots, -3.3, \dots, 0, \dots, 1, \dots, 100.23, \dots$, i.e., the whole numbers, and the decimals between them. Catherine claims the order on $\mathbb R$ is dense, which is to say, if we choose any two distinct numbers in the real numbers, we can find a distinct number between them.

Catherine would construct her proof in the manner laid out in the introductory monograph \emph{When is a Proof?}~\cite{davey_when_2009}, in Table \ref{tab:steps}, provided to undergraduate mathematics majors at La Trobe University. These steps comprise \textbf{formal} and \emph{informal} mathematical work, showing that mathematical \emph{work} comprises more than the \emph{advertising}, as it is labelled in the reproducibility spectrum presented in Figure \ref{fig:repro-spectrum}. In the case of pure mathematics, the advertising would be the paper that outlines the proof, the formal mathematical argument, but the informal work is left out.

\begin{table}
\caption{The steps in the making of a proof from Brian A. Davey's primer, \emph{When is a Proof?}~\cite{davey_when_2009}. The formal steps that contribute to the final proof are in \textbf{bold}, the hidden informal work, in \emph{italics}. These steps are summarised in terms of $p \implies q$ in Table~\ref{tab:pq-steps}.} \label{tab:steps}
    \centering
\footnotesize{    
    \begin{tabular}{c l p{10cm}l}
    \toprule
  Step -1 & $\qquad$&  \textbf{Translate the statement to be proved into ordinary English and look up appropriate definitions.}\\
  Step 0  && \textbf{Write down what you are asked to prove. Where appropriate, isolate the assumptions, p, and the conclusion, q.}\\
  Step 1 &&  \textbf{Write down the assumptions, p: “Let . . . . . . ”}\\
  Step 2 &&  \textbf{Expand Step 1 by writing out definitions: “i.e. . . . . . . ”}\\
  \midrule
  Step 3 &&  \emph{Write down the conclusion, q, which is to be proved: “To prove: . . . . . . ”}\\
  Step 4 && \emph{Expand Step 3 by writing out definitions: “i.e. . . . . . . ”}\\
  Step 5 && \emph{Use your head: do some algebraic manipulations, draw a diagram, try to find the relationship between the assumptions and the conclusion.}\\
  \midrule
  Step 6 && \textbf{Rewrite your exploration from Steps 3, 4 and 5 into a proof. Justify each statement in your proof.}\\
  Step 7 && \textbf{The last line of the proof: “Hence q.”}\\
  \bottomrule
  \end{tabular}
 } 
\end{table}

Catherine presents the following proof to Hal to show the order on real numbers, $\mathbb R$, is dense.   

\begin{proof}{The order on $\mathbb R$ is dense.} 
Let $x < y$ in $\mathbb R$. Let\footnote{
In mathematics, we read $:=$ as `be defined as', $\implies$ as `implies', and $<$ as `less than but not equal to'.
} $z := \frac{x + y}{2}$. To see that $x < z < y$, we begin with $x < y$, so, $x + x < x + y$ and $x +  y < y + y$, which gives,

$$
\begin{array}{ccccccc}
\centering
    &x  + x & < &x + y & < &y + y\\
    &&&&&\\
    \implies &\displaystyle\frac{x + x}{2} & < &\displaystyle\frac{x + y}{2} & < &\displaystyle\frac{y + y}{2}\\
    &&&&&\\
    \implies &\displaystyle\frac{2x}{2} & < &\displaystyle\frac{x + y}{2} & < &\displaystyle\frac{2y}{2}\\
        &&&&&\\
     \implies & x & < & \displaystyle\frac{x + y}{2} & < & y\\
         &&&&&\\
     \implies & x & < & z & < & y,
\end{array}
$$

since $z = \frac{x + y}{2}$, as required.\qed 
\end{proof}

\begin{figure}
    \centering
    \includegraphics[width=0.46\textwidth]{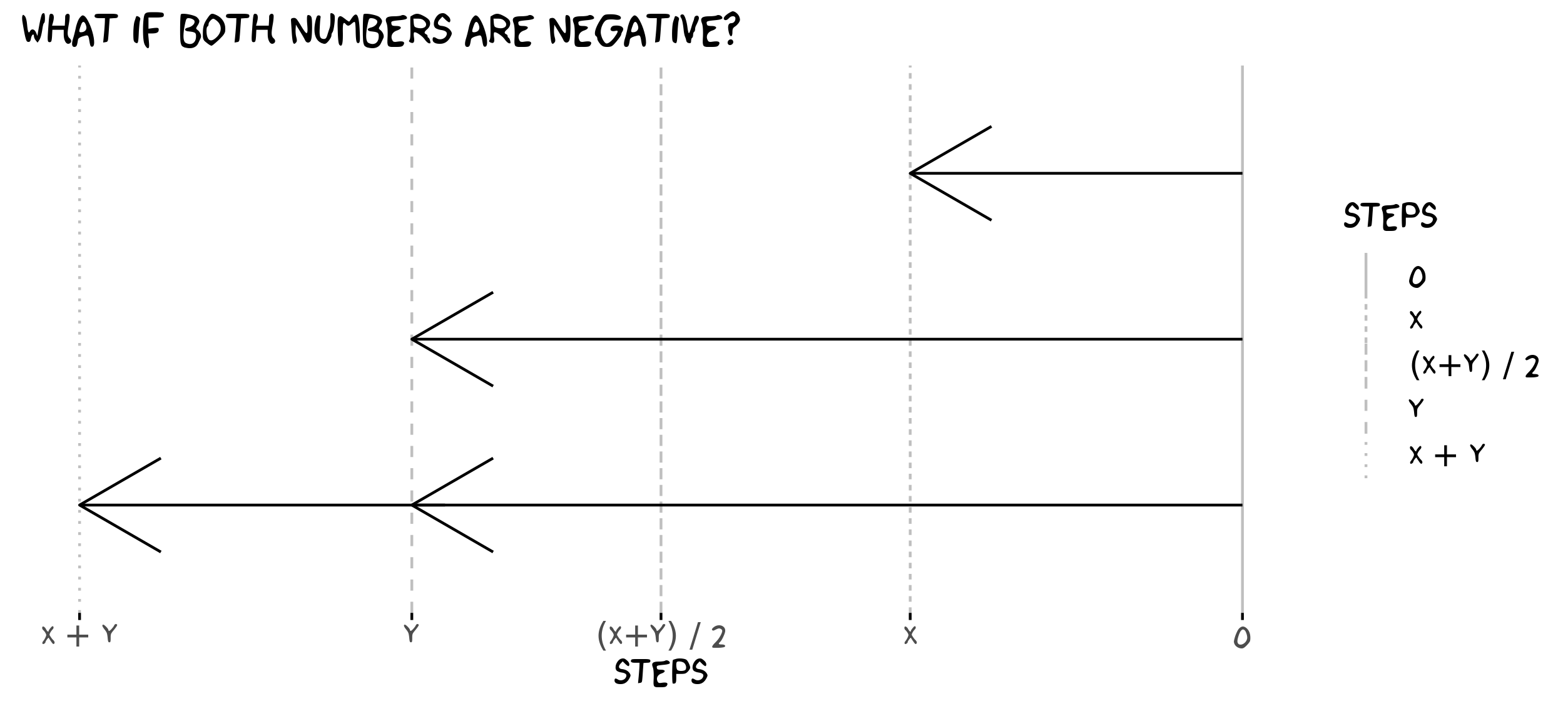}
    \includegraphics[width=0.46\textwidth]{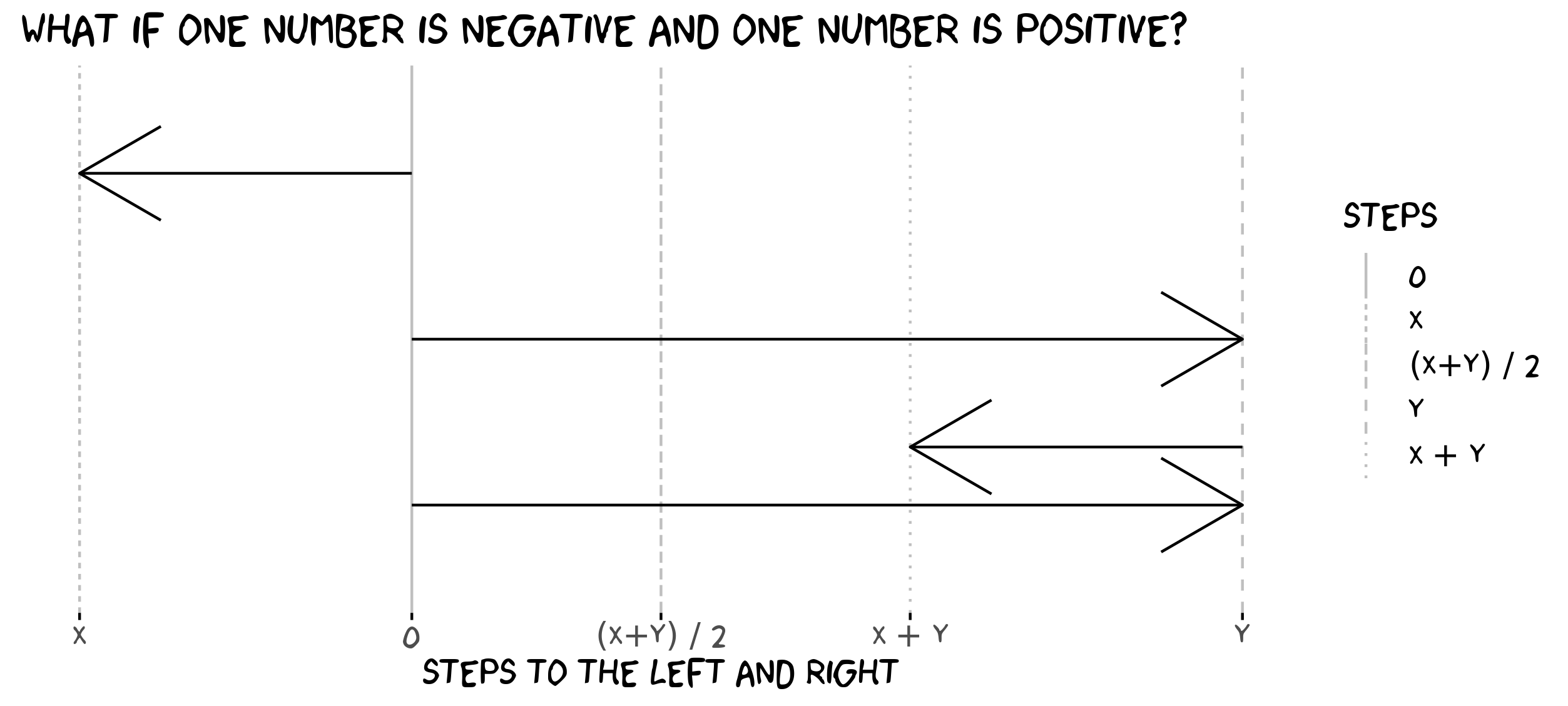}
     \caption{On the left, Hal might begin to verify his understanding of $+$ by first considering the case where both numbers are negative, $x, y < 0$. In this case, we might think of $+$ as combining $x$ steps to the left with $y$ steps to the left. The halfway point $\frac{x + y}{2}$, falls in the middle of the two arrows laid side by side. On the right, Hal considers the case where $x < 0$, $y > 0$ and $|x| < |y|$. Here $x + y$ can be thought of as $y$ steps to the right and then $x$ steps to the left. Again, the halfway point $\frac{x + y}{2}$ falls halfway between 0 and $x + y$.}\label{fig:neg-steps}

 \end{figure}

Catherine presents the formal proof, the science that in Figure \ref{fig:repro-spectrum} is described as the advertising, a subcomponent, of the compendium of research she created in order to arrive at this argument. Hal wishes to verify the results and investigate whether Catherine merely reproduced her father's reasoning. In the case of proof, what is published is the formal argument, but as the steps in Table \ref{tab:steps}, this is not all of what makes a proof. We could think of the steps presented in Table \ref{tab:steps} in terms of a mathematical statement $p \implies q$, which we read as $p$ \emph{implies} $q$. In Table \ref{tab:pq-steps}, we summarise the \emph{Steps in the making of a proof} in terms of $p \implies q$.
\begin{wraptable}{l}{0.3\textwidth}
    \caption{Steps in the making of a proof, shown in Table \ref{tab:steps}, in terms of $p \implies q$.}
    \label{tab:pq-steps}    \centering
\footnotesize{    
    \begin{tabular}{ccl}
    \toprule
    0   &$\qquad$&  $p \implies q$  \\
    1   && \textbf{Assume $p$.}\\
    2   && \textbf{Define $p$.}\\
    \midrule
    3   && \emph{State $q$.}\\
    4   && \emph{Define $q$.}\\
    5   && \emph{Work.}\\
    \midrule
    6  && \textbf{Formalise work.}\\
    7  && \textbf{So, $q$.}\\
    \bottomrule
    \end{tabular}
} 
\end{wraptable}
 We begin, step 0; we state what we wish to prove $p \implies q$ in plain English. We wish to show the real numbers, $\mathbb R$, are dense; i.e., for all $x < y$ in $\mathbb R$, there exists $z$ such that $x < z < y$.\\ 
 Step 1, we \textbf{assume} $p$ is true. We assume we have two distinct numbers $x$ and $y$ in $\mathbb R$ with $x < y$; i.e., $x$ is less than $y$, and $x$ is not equal to $y$. Step 2, nothing to define as we are familiar with $<$ and $\mathbb R$. \\
 Step 3, we wish to \emph{prove} $q$; Step 4, i.e., we need to show there exists $z$ in $\mathbb R$ such that $x < z < y$. Now, Catherine has offered a solution $z := \frac{x + y}{2}$ that Hal wishes to verify.\\ 
 Step 5, Suppose Hal asks, what if both $x$ and $y$ are negative numbers? Is it still true that $x < z < y$? Hal might verify his understanding of $+$ by thinking about positive and negative numbers as steps taken to the left or the right. In Figure \ref{fig:neg-steps}, Hal considers the case where both numbers are negative, $x, y < 0$. In this case, we have $x$ steps to left, and $y$ steps to the left, which we imagine as arrows of appropriate length. If we lay both arrows end to end, we see the number of combined steps to the left. If we consider the half-way point of $x$ and $y$ laid beside each other, $\frac{x + y}{2}$, we see this falls between where the arrow heads of $x$ and $y$ fall. 

Now Hal can flip the arrows in the opposite directions to construct an argument for if both numbers were positive, $x, y > 0$. 

But then Hal asks in Figure \ref{fig:neg-steps}, what if one number were positive and one number were negative? Is $\frac{x + y}{2}$ still halfway between? Let us assume, as mathematicians say, without loss of generality that the magnitude of $x$ is strictly less than $y$, that is $|x| < |y|$, the number of steps in $x$ is less than the number of steps of $y$. Hal now considered where one would end up if one took $y$ steps to the right and then $x$ steps to left. He checks that he does not need to consider two cases, as he would end up in the same place if he took $x$ steps to the left and then $y$ steps to the right. Again, $\frac{x + y}{2}$ falls between where he would start and where he would end. 

Now Hal has verified his understanding of $+$, which may or may not be the way that Catherine arrived at her result, but after this work he is capable of fully reproducing the mathematical result presented. He reads the proof Catherine has provided, and verifies Steps 6, and Step 7. Catherine has proved that the order on $\mathbb R$ is dense. With this proof, as with the proof presented in Section \ref{sec:prove-it}, Hal cannot disqualify the possibility that Catherine merely reproduced her father's work. 

Even in these relatively simple proofs, Step 5, the informal work of verification and understanding vastly outweighs what goes into the formal proof.  But these toy examples belie a process of redefinition and re-examination, as illustrated in the discussion within a hypothetical mathematics classroom that forms the narrative of Lakatos' \emph{Proofs and Refutations}~\cite{lakatos_proofs_2015}. We now move to a recently published proof to illustrate this process of redefinition.

\subsubsection{In the combat conditions of new mathematics}\label{sec:combat}

Suppose, now, that Catherine's proof were for the theorem pertaining to quasi-primal algebras, presented in the recent publication `The homomorphism lattice induced by a finite algebra'~\cite{davey_homomorphism_2018} in  \emph{Order}, a mathematics journal devoted to `original research on the theory and application of ordered sets'. In addition to the informal work demonstrated by the proof that the order on  $\mathbb R$ is dense, the making of this proof involved a redefinition of the result proved, through a process writing several proofs. In terms of Table \ref{tab:pq-steps}, initially a result was considered, $p \implies q$. A proof was written for this result. At this point the mathematicians realised, however, that the converse could be shown, that is, $q \implies p$. And so, a proof was generated for a new result, $p \iff q$. In the case of this proof, the act of writing the proof itself redefined the result in question. In the combat conditions of new mathematics, the process of writing a proof is doing mathematical science, and involves a great deal more work than is presented in the advertising of the science. 

Hal may require graduate-level knowledge of abstract algebra to reproduce this proof, but as a professional mathematician, this is not a great leap. More challenging the proof may be, but the process of reproduction would be similar. Even if this were the proof, Hal would not know if Catherine merely reproduced, as he did, her father's proof. 

But what if Catherine were posing her mathematical question computationally? Would Hal be able to reproduce her results? 

\subsection{Is computational mathematics mired in proof methodology?}

When we are exploring and answering mathematical questions in a computational environment, we consider some aspects of our work to be \textbf{formal} and some \emph{informal}. But in omitting the greyed \emph{informal} work in Figure \ref{fig:repro-spectrum}, are we still approaching compendia of research from the perspective of a blackboard mathematician?  

Given we use statistics in most science, arguably most scientific questions are posed, to some extent, mathematically. The output format, a published paper, remains similar in format to mathematics of the pre-computer age. But the informal work of answering  mathematical questions has changed significantly. Now that much work is done computationally, there are multiple research outputs that comprise the compendium of science that produces the published paper.

Let us now suppose that Catherine had a statistical estimator for a population parameter of interest. That is, Catherine has an equation that, given some data, she can approximate some value about the population, such as an overall average. Let us further suppose, as is increasingly common, that she does not have a closed-form solution, meaning she cannot write out a mathematical argument in the traditional sense. Instead, she demonstrates the estimator's performance through simulation studies. 

Now suppose Hal challenges Catherine to prove that she created the science that produced the paper. Given what is on the piece of paper, how can Hal know that Catherine's code does what she said it does? It is unclear what assumptions were made, about, say, sample size and distribution. How can Hal verify her results? Through adopting research software engineering principles, Catherine can facilitate a process akin to proofs and refutations, the redefinition described in the Section \ref{sec:combat}, The combat conditions of new mathematics. The process of redefinition is transcribed by version control, but further to this, the software itself provides a modular framework, such as a theorem in mathematics, for future work to scaffold and extend. New software can be developed that either extends, or redefines the existing software. One analogous way this is occurring is in the rise of metapackages, such as \verb|tidyverse::|~\cite{wickham_tidyverse_2017} and \verb|metaverse::|~\cite{metaverse_2019}, that collect software to solve particular problems in an opinionated~\cite{parker_opinionated_2017} manner, that guide the end-user to what the creators consider to be good enough practice. This is analogous to classes of mathematics, such as group theory or analysis, that collect results, theorems, that rely upon each other, and where certain underlying assumptions, such as the \emph{Axiom of Choice}\footnote{
Turning to the bible of algebra, \emph{Lattices and Order}~\cite{davey_introduction_2002-1}, we learn the \emph{Axiom of Choice} `asserts that it is possible to find a map which picks one element from each member of a family of non-empty sets'. 
}, are made.  

How are contemporary researchers answering mathematical questions? Alex Hayes, current maintainer and one of the many authors of \verb|broom::|~\cite{robinson_broom_2019}, an open source R package that amalgamates hundreds of contributions towards providing a suite of tools that  \emph{tidily}\footnote{From Wickham's \emph{Tidy data}~\cite{hadley_wickham_tidy_2014}, we describe data as \emph{tidy} if 
\begin{enumerate}
    \item Each variable forms a column.
    \item Each observation forms a row.
    \item Each type of observational unit forms a table. 
\end{enumerate}
}~\cite{hadley_wickham_tidy_2014} extract statistical model information from R algorithms, recently noted the underdeveloped nature of the implementation of statistical algorithms~\cite{hayes_testing_2019}: 

\begin{quote}
    `In practice, most people end up writing a reference implementation and checking that the reference implementation closely matches the pseudocode of their algorithm. Then they declare this implementation correct. How trustworthy this approach is depends on the clarity of the connection between the algorithm pseudocode and the reference implementation.' 
\end{quote}

This is not to carp upon diligent scientists; we need to do far more to support the software engineering principles we expect from mathematical scientists~\cite{nowogrodzki_how_2019}. Mathematicians are trained to provide enough work such that the hidden steps illustrated in italics in Table \ref{tab:steps} can be reproduced by their target audience. The detail of mathematical work shown is tempered for level of the audience, but the same process described in bold in Table \ref{tab:steps} is the same. But, does the workflow Alex describes above equip the target audience with enough information such that they can understand all the details of the entire argument put forward? 

Code has the appearance of being highly logical, it's easy to assume it's infallible; and whilst the logic of the code is robust, the pipeline that carries the algorithm to implementation may be susceptible to compromising factors, with typos being just one example of inadvertent error. 

Because code appears so logical, we assume it is analogous to proof for our intended audience to follow. But we were trained to leave out the informal messy thinking work associated with mathematics; trusting the formal argument provides enough information to verify and reproduce the mathematics. Does our code do what we think it does? In addition to providing the  research outputs in the spectrum of reproducibility, Figure \ref{fig:repro-spectrum}, we posit mathematical science should adopt the software development practice of unit testing, to ensure the mathematical results can be verified and reproduced.    

\section{Testing}

\begin{wraptable}{}{0.5\textwidth}
    \caption{Percentage of R packages in repositories that have unit tests included. These results are from Jim Hester's \href{https://www.rstudio.com/resources/webinars/covr-bringing-test-coverage-to-r}{\underline{presentation}} on covr:: in September 2016 ~\cite{hester_covr_2016}.}
    \label{tab:hester}
    \footnotesize{
    \begin{tabular}{rlllr}
        \toprule
        Repository $\quad$ & $\quad$ & Tests $\quad$ & Total $\quad$ &    \\
        \midrule
        CRAN & &  2091 &  9772 & 21\%\\
        Bioconductor & & 449 & 1258 & 36\%\\
        rOpenSci & & 84 & 146 & 58\%\\
        \midrule
         & & 2624 & 11,176 & 24\%\\
         \bottomrule
    \end{tabular}
    } 
\end{wraptable}
Testing is the software engineering tool that is provides a key piece of the correspondence between scientific claim and programming. Just as the Curry-Howard isomorphism expresses proofs-as-programs to link mathematics and programming, we argue that tests are are the link between scientific claims more generally and programming. In a test the researcher isolates a scientifically meaningful part of their code, and creates a witness so that others can easily see that the code does what the researcher intends it to do. In this section we consider a `vital'~\cite{wickham_r_2015} research output, testing, that it is unlikely the mathematical scientist has been trained in. There are many such under-formalised skills represented in Figure \ref{fig:repro-spectrum}\footnote{Indeed, the natural consequence of questioning how we practice mathematical science is how we train the next generation of practitioners. Important, however this may be, this is beyond the scope of this manuscript.}. In 2016, a quarter of packages on R package archives CRAN, Bioconducter, and rOpenSci, included tests, a repository by repository breakdown of this is shown in Table \ref{tab:hester}.

Now, Hayes advises people against using untested software~\cite{hayes_testing_2019}. It is alarming that, by this logic, we would be \textbf{insane} to use \emph{three quarters} of packages available. But Hayes continues, `You have two jobs. The first job is to write correct code. The second job is to convince users that you have written correct code'~\cite{hayes_testing_2019}. The disconnect here suggests a failure to communicate broadly the importance of testing of algorithms in the dissemination of research. As researchers, we believe our science is as reproducible as a traditional mathematical proof; however, the growing literature of the replication crisis demonstrates we have not succeeded in rendering our science reproducible.

rOpenSci's review system recommends using the \verb|covr::|~\cite{hester_covr_2018} package to measure how the code behaves with different expected outputs. From the creator of \verb|covr::|, we obtain the following definition of test coverage.

\begin{quote}
`Test coverage is the proportion of the source code that is executed when running these tests'~\cite{hester_covr_2018}.
\end{quote}

\subsection{What is a test?}

Tests demonstrations that a given input produces an expected output. They are grouped contextually in a file; the context being a certain aspect of the algorithm that should be tested ~\cite{wickham_r_2015}. An example of a context for a test is the question, does a given function return the expected result for different inputs? Each test comprises a collection of expectations. Each expectation runs a function or functions from the package, and checks the returned output is as expected. In this case, we have a test for the \verb|expect_equal| function: one expectation checks the  function successfully runs when given equal inputs, and another expectation checks that the function fails when passed two non-equal inputs.

An example test from the \verb|testthat::|~\cite{wickham_testthat_2011} contains two expectations. 
\begin{verbatim}
    test_that("basically principles of equality hold", {
  expect_success(expect_equal(1, 1))
  expect_failure(expect_equal(1, 2))
})
\end{verbatim}


\subsection{How good are we at \emph{good enough} testing?}

A response to the replication crisis has been to examine \emph{questionable research practices}~\cite{fraser_questionable_2018}, frequently borne of tradition and convention within different disciplines, deviate from evidence-based best-practice research methodology. We suggest it is a questionable research practice to draw conclusions about the efficacy of statistical estimators from untested code. 

Given only a quarter of R packages have unit tests associated with them, we are falling short of best practice in scientific computing~\cite{wilson_best_2014}. In a recent assessment of what constitutes \emph{good enough} practice in scientific computing~\cite{wilson_good_2017}, unit testing was not included. However, for mathematical science, where the algorithms implemented and the code written is often complex, we suggest that unit testing should be considered good enough practice, in spite of the additional learning curve. With the backdrop of the replication crisis, it is crucial we have confidence in the algorithms we implement.  

\subsection{Analysis of testing code in R packages}

\begin{figure}
    \centering
    \includegraphics[width = \textwidth]{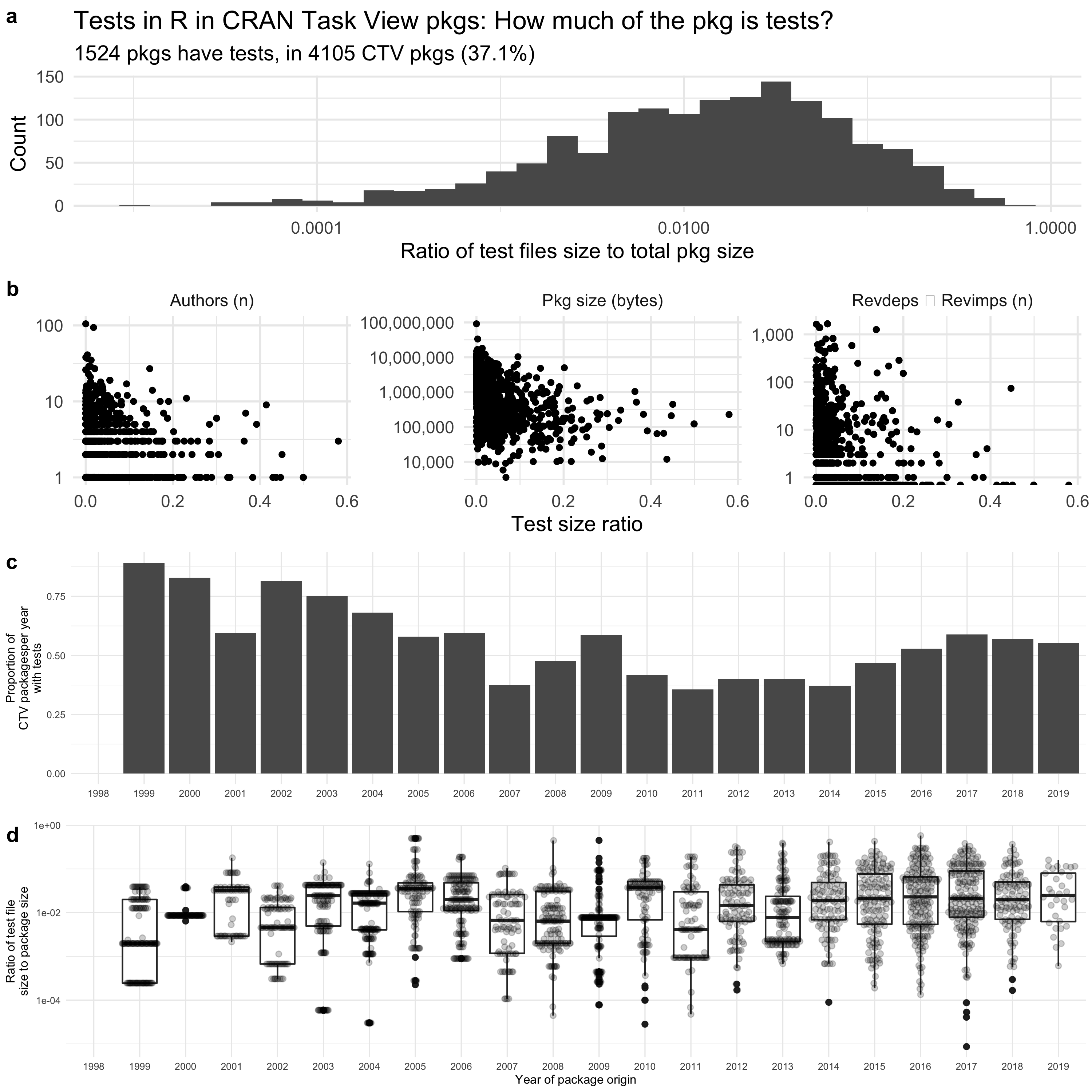}
    \caption{This panel shows some basic details of tests in R packages listed in CRAN task views~\cite{zeileis_cran_2005-1}. The measure of interest, \emph{test size ratio}, was calculated by dividing the test file size with the overall package source file size from the unofficial CRAN mirror on GitHub. This is a rough indicator of test coverage, future work should consider more precise metrics such as those produced by the covr:: package. a) the distribution of the ratio of test file size to total package size, test size ratio. b) scatter plots demonstrate the relationship between test size ratio and number of authors, overall package size, and number of packages imported and calling the package, respectively. c) the proportion of all task view packages that contain tests over time. d) boxplot detailing the distribution of file size ratio over time.}
    \label{fig:test-analysis}
\end{figure}

So, what packages have tests? We provide a preliminary analysis of tests in CRAN packages in Figure \ref{fig:test-analysis}. The code and data used to generate the results presented here are openly available at \href{https://github.com/softloud/proof}{https://github.com/softloud/proof}.

We provide analysis for packages associated with CRAN task view~\cite{zeileis_cran_2005-1}, opinionated~\cite{parker_opinionated_2017} collections of R packages that are relevant to a particular type of statistical analysis, maintained voluntarily by experts in their respective fields~\cite{zeileis_cran_2005-1}. CRAN task views provide a convenient taxonomy of R packages for a preliminary exploratory analysis of patterns of test use among R package authors.

Packages listed in a task view are may be interpreted by users as more stable and trustworthy than other packages, because they have passed some kind of inspection by  maintainer of the task view who listed the package (however the review and curation process is not open or documented). And yet, even amongst the 4105 packages associated with task views, 1524 packages were without tests; 37 per cent of packages associated with CRAN task view were without tests. 

The proportion of task view packages with tests has fallen over the last decade. This does not seem surprising given the uptake of R amongst communities of researchers in applied sciences with little formal programming and computer science training, such as psychology and ecology.

\begin{figure}
    \centering
    \includegraphics[width = \textwidth]{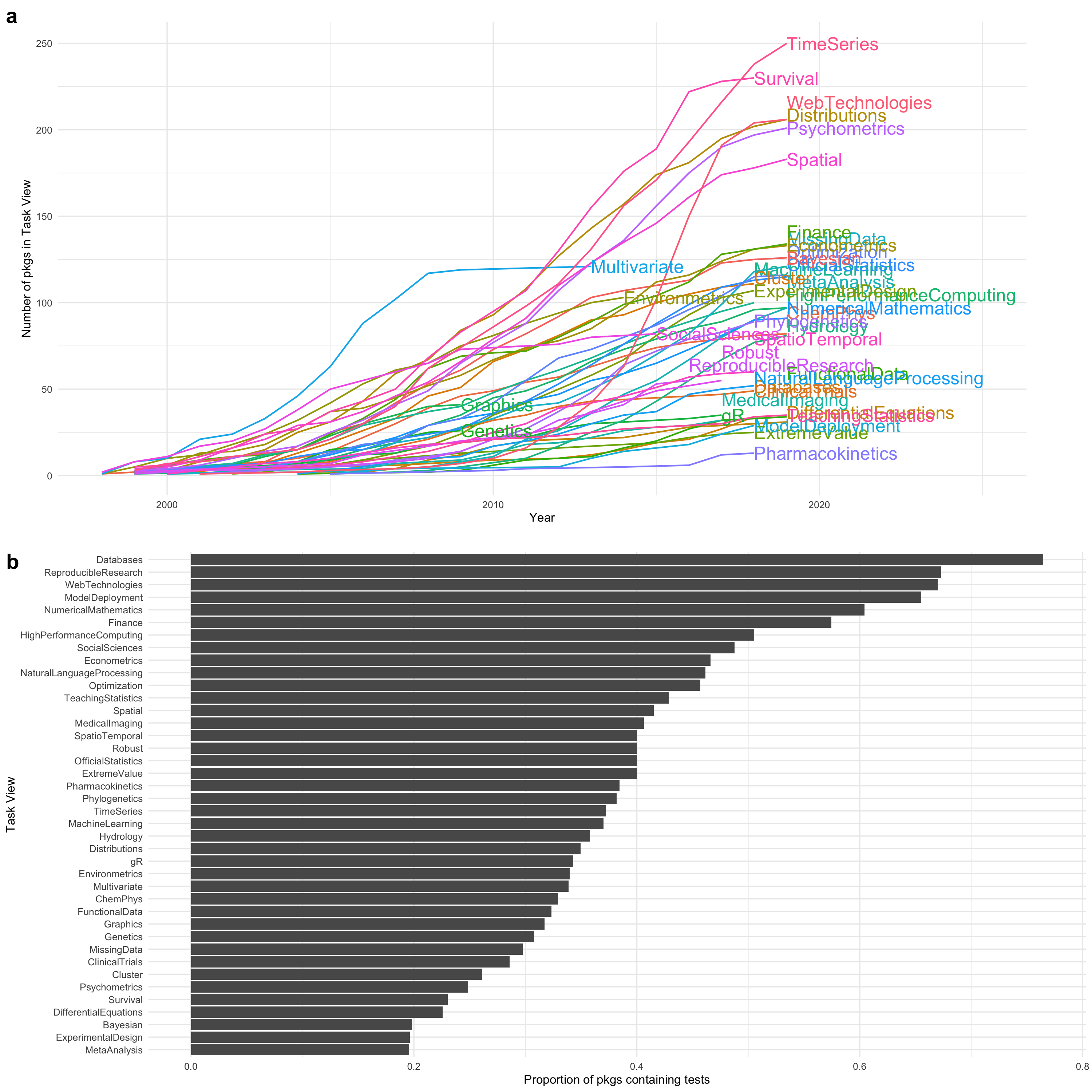}
    \caption{a) shows the change in the number of packages in each CRAN task view over time. b) shows the proportion of packages in each CRAN task view that have tests.}
    \label{fig:ctv-prop}
\end{figure}

Figure \ref{fig:ctv-prop} shows that there is wide variation in test coverage. Even the largest and fastest growing CRAN task views have very different proportions of packages with tests (Survival, about 0.23, compared to Web Technologies about 0.66). We find few clear patterns in the presence of tests over time, between different CRAN task views, and with metadata such as the number of authors, the size of the package and the centrality of the package (as measured by the union of the number of reverse dependencies and reverse imports). Based on these data, we suggest there is much work to be done in developing methods and opinionated tools that guide users towards good enough practices.

\section{Tempered uncertainty and computational proof}

It's easy to lie with statistics, but it's even easier without them~\cite{murray_how_2005}. In a computational experimental setting, we often cannot achieve the satisfying precision offered by a proof. We can, however, adopt good enough practices in sharing and testing code to increase confidence in our scientific conclusions. It may not be possible to provide the rigour of a closed-form mathematical solution, but we can aim to temper the uncertainty, and bolster confidence, in computational arguments via automated testing, version control, and other computational outputs. 

We suggest there is much work to be done in developing good enough practices~\cite{wilson_good_2017} we can ask mathematical scientists to adopt. Indeed, the question of good enough practice can be posed for each research output. Less than offering answers, this manuscript seeks more to suggest there is a rich line of inquiry~\cite{nowogrodzki_how_2019} in the relationship between scientific truth, mathematical proof, and computational reproducibility and rigour. 

\subsection{Coda}

Returning to Catherine and Hal from Auburn's \emph{Proof}~\cite{auburn_proof_2001}, we can now imagine her as computational mathematician who provides a compendium of reproducible research. To demonstrate the rigour of her computational work, she would provide unit tests for the algorithms she had implemented. Catherine would share her work openly via her GitHub or similar repository, where the development of her ideas would be timestamped and recorded. The structure of her research compendium of would be automatically standardised via a tool such as \verb|rrtools::|~\cite{marwick_rrtools_2018}. At publication, her compendium would be deposited on a trustworthy, DOI-issuing repository for others to link to and cite. 

And she would feel safe asking questions about good enough practice~\cite{wilson_good_2017}, and how to avoid questionable research practices~\cite{fraser_questionable_2018}, because there is an understanding in the community that no one is trained in all these things, so we are all always learning.  

There would be no struggle, as there was in Auburn's play, to show that the mathematician who created these research outputs was Catherine. But that wouldn't matter - she and Hal would be having far too much fun collaborating on the next question. 


\bibliographystyle{splncs04}
\bibliography{references.bib}

\end{document}